\documentclass[12pt]{amsart}
\usepackage[T1]{fontenc}
\usepackage[utf8]{inputenc}
\usepackage{amssymb}
\usepackage{amsthm}
\usepackage{float}
\usepackage{longtable}
\usepackage{amsmath}
\usepackage{mathrsfs}
\usepackage{nicematrix}
\usepackage{yfonts}
\usepackage{mathtools}
\usepackage{tikz-cd}
\usepackage{xcolor}
\usepackage{fullpage}
\usepackage{graphicx}
\usepackage[parfill]{parskip}
\usepackage{hyperref}
\usepackage{biblatex}
\hypersetup{
  colorlinks   = true,    
  urlcolor     = black,    
  linkcolor    = black,    
  citecolor    = black      
}
\newtheorem*{theorem*}{Theorem}
\newtheorem{theorem}{Theorem}

\newtheorem{proposition}[theorem]{Proposition}

\newtheorem{definition}[theorem]{Definition}

\numberwithin{theorem}{section}

\newcommand\rr{\rightarrow}
\newcommand{\PP}{\mathbb P}

\newcommand{\QQ}{\mathbb{Q}}

\newcommand{\CC}{\mathbb{C}}

\newcommand{\mc}{\mathcal}

\newcommand{\wdt}{\widetilde}

\begin{document}
    
\title{Degeneration of Calabi-Yau threefolds with Galois action induced by modular form of weight 3}
\author[]{Marcin Oczko}

\begin{abstract}
We study a one-parameter family of double octic Calabi-Yau threefolds with two degenerations of type $K$. For each of these degenerations we construct a semistable model and compute the limiting mixed Hodge structure. We show that the $
\ell$-adic Galois representation on the limiting cohomology is isomorphic (up to semisimplification) to the direct sum of the Galois representation attached to a weight-three modular form associated with a singular K3 surface and its Tate twist. This gives an explicit realization of the expected relation between type $K$ degenerations of Calabi-Yau threefolds and K3 surfaces.

\end{abstract}

\thanks{The author was supported by the National Science Center, Poland, grant no. 2025/57/N/ST1/05021.}

\maketitle


\section{Introduction}

In this paper we study a one-parameter family of Calabi-Yau threefolds with two conifold points and two points of type $K$. This family has the following striking property: point counts of the singular fibers of type $K$ over finite fields are related by explicit formulas to the Fourier coefficients of weight-three modular forms attached to singular K3 surfaces. This suggests that part of the limiting mixed Hodge structure of the degenerate fibers should arise from the transcendental cohomology of K3 surfaces. We prove that this is indeed the case. Using the modular form notation from LMFDB \cite{lmfdb}, the main result may be stated as follows.

\begin{theorem*}[Theorem \ref{main thm}]
    Let $\mc{Y}$ be a family of Calabi-Yau threefolds (family No. 70 from \cite{CYK}), with two singular fibers $Y_0,Y_\infty$ of type K. Let $f_1,f_2$ be two modular forms 32.3.d.a and 16.3.c.a. of weight-three associated to K3 surfaces \cite{lmfdb}. The Galois representation on the limiting mixed Hodge structure of $Y_0$ (resp. $Y_\infty$) is isomorphic up to semisimplification to the direct sum of the Galois representation associated to $f_1$ (resp. $f_2$) and its Tate twist.
\end{theorem*}

The proof is based on an explicit semistable model of the family. Using the monodromy weight spectral sequence, we compute the limiting mixed Hodge structure and show that its two-dimensional graded pieces are induced by the transcendental lattice of the corresponding singular K3 surface.

Degenerations we study are locally over a complex disk $\mc{X} \rr \Delta$ with a singular fiber $X_0$. The monodromy operator $T \colon H^d(X_t,\QQ_\ell) \rr H^d(X_t,\QQ_\ell)$ is given by parallel transport around 0. By a fundamental result \cite{SGA7} the operator $T$ is \emph{quasi-unipotent}, i.e. $(T^k-I)^n=0$ for some $k,n$. Moreover, for \emph{semistable} degenerations this operator is unipotent \cite{Landman}. From studying the monodromy operator $T$ we can infer certain properties of the degenerate fiber. In our case we will use it to study families of Calabi-Yau varieties.

For two-dimensional Calabi-Yau varieties, i.e. K3 surfaces, we distinguish three types of degenerations depending on the unipotency index of $T$. By the theorem of Kulikov \cite{Kul77,Kul81} and Persson and Pinkham \cite{PP81}, the unipotency index determines the type of the degeneration. In particular a family of K3 surfaces has a smooth filling if the monodromy action is trivial.

In the arithmetic context Chiarellotto, Lazda and Liedtke proved a generalization of the Néron–Ogg–Shafarevich criterion for K3 surfaces \cite{CLL}, which is an arithmetic analogue of the Kulikov-Persson-Pinkham theorem. Their result gives conditions for the Galois representation associated with a K3 surface in characteristic zero, which determine if it has a good reduction modulo a prime or what kind of singularities might appear.

For three-dimensional Calabi-Yau manifolds analogues of both classical and arithmetical characterizations of singularities of degenerate fibers no longer hold. There is a family of Calabi-Yau threefolds in characteristic zero constructed by Cynk and van Straten \cite{CvS21}, which has a trivial monodromy, but no smooth filling. There is also a Calabi-Yau threefold in characteristic zero with unramified Galois representation yet admitting no good reduction to characteristic $p$ \cite{CO24}.

For families of Calabi-Yau threefolds with $h^{1,2}=1$, we have the following possibilities for the Jordan form of the local monodromy action.

\begin{table}[H]
\begin{tabular}{c | c}
type of singularity & Jordan form of M \\ \hline
F & $\begin{pmatrix}
        1 & 0 & 0 & 0\\
        0 & \zeta_{n_1} & 0 & 0 \\
        0 & 0 & \zeta_{n_2} & 0\\
        0 & 0 & 0 & \zeta_{n_1} \zeta_{n_2}
    \end{pmatrix}$  \\ \hline
C & $\begin{pmatrix}
        1 & 0 & 0 & 0\\
        0 & \zeta_{n} & 1 & 0 \\
        0 & 0 & \zeta_{n} & 0\\
        0 & 0 & 0 & \zeta_{n}^2
    \end{pmatrix}$  \\ \hline
MUM & $\begin{pmatrix}
        1 & 1 & 0 & 0\\
        0 & 1 & 1 & 0 \\
        0 & 0 & 1 & 1\\
        0 & 0 & 0 & 1
    \end{pmatrix}$ \\ \hline
K & $\begin{pmatrix}
        1 & 1 & 0 & 0\\
        0 & 1 & 0 & 0 \\
        0 & 0 & \zeta_{n} & 1\\
        0 & 0 & 0 & \zeta_{n}
    \end{pmatrix}$    
\end{tabular}
\end{table}

Here $\zeta_n$ denotes the primitive $n$-th root of unity.

It is expected that at MUM points we have a maximal degeneration, at C points a (singular) rigid Calabi-Yau threefold and at F points (after a finite cover) a smooth manifold. 

The case of degenerations of type K is the least studied. The name itself comes from the expected link with K3 surfaces as explained in \cite{CvS21}.

\section{Families of double octics}

Throughout this paper we work with a class of Calabi-Yau threefolds called \emph{double octics}. They are constructed as a crepant resolution of a double covering $X$ of $\PP^3$ branched along a union of 8 planes $S=S_1 \cup ... \cup S_8$. Such a crepant resolution exists if no six planes in $S$ intersect and no four planes contain a common line.

Singularities of $X$ are in one-to-one correspondence with intersections of components in its branch locus. We identify the following types of singularities in the branch locus of $X$.

\begin{itemize}
    \item Lines contained in $k$ different planes, denoted by $l_k$. We allow lines of types $l_2$ and $l_3$.
    \item Isolated $j$-fold points contained in exactly $k$ triple lines denoted as $p^k_j$. We allow isolated singular points of types $p^0_3,p^0_4,p^1_4,p^0_5,p^1_5$ and $p^2_5$.
\end{itemize}

Singularities of a double cover can be resolved by a sequence of blow-ups of components in the singular locus in the following order:

\begin{enumerate}
    \item blow-up of  all fivefold points,
    \item blow-up of all triple lines,
    \item blow-up of all fourfold points,
    \item blow-up of all double lines (in some chosen order).
\end{enumerate}

This procedure yields a sequence $X_i \rr P_i$, $i=0,...,n$ of double coverings with branching loci $S^{(i)} \subset P_i$ and blow-ups $\sigma_i \colon P_i \rr P_{i-1}$, $i=1,..n$ along smooth centers $C_i$ with $P_0 = \PP^3$, $S^{(0)}=S$, $X_0=X$ and $Y:=X_n$ smooth. If the multiplicity of $S^{(i)}$ at a generic point of $C_i$ is odd (triple lines and quintuple points) then $S^{(i)}=\wdt{S^{(i-1)}}+Exc(\sigma_i)$, if the multiplicity of $S^{(i)}$ at a generic point of $C_i$ is even (double lines and quadruple points) then $S^{(i)} = \wdt{S^{(i-1)}}$.

Consider the following family of double covers of $\PP^3$ (Arr. No. 70 in \cite{CYK}) with the parameter $w = (w_0:w_1)\in \PP^1$.

\begin{multline}\label{eq1}
\mc{X}:=\{((x:y:z:t:u),(w_0:w_1))\in \PP(1,1,1,1,4)\times\PP^1: \\
xyzt(x - y + z) (y - z - t) (x - y - t) (w_1x + w_0y)=u^2\}
\end{multline}

We denote intersections $S_{i_1 \cdots i_n}= S_{i_1}\cap \cdots \cap S_{i_n}$. For $w \neq 0,-1,-2, \infty$ the fiber $X_w \in \mc{X}$ has the following singularities:

\begin{itemize}
    \item $p^0_4$: $S_{1456},S_{2346},S_{2567},S_{3457}$,
    \item $p^1_4$: $S_{1268}$,
    \item $p^1_5$: $S_{12358},S_{12478}$,
    \item $l_3$: $S_{128}$.
\end{itemize}

For $w=-1,-2$ we get octic arrangements with additional singularities (Arr. Nos. 3 and 69 in \cite{CYK}). For $w=0$ (resp. $w=\infty$) planes $S_1$ and $S_8$ (resp. $S_2$ and $S_8$) coincide, hence these fibers cannot be resolved to a double octic. 

We define a family $\mc{Y}$ as a completion of a fiberwise resolution of generic fibers of $\mc{X}$. Fibers $Y_w$ of $\mc{Y}$ for $w \neq 0,-1,-2,\infty$ are smooth Calabi-Yau threefolds. The behavior of $\mc{Y}$ near a singular point is dictated by the local monodromy described by Riemann symbols (local exponents) of the Picard-Fuchs operator (computed in \cite{Orph}, Section 5) which in this case are the following.

\begin{table}[H]
\begin{tabular}{cccc}
0 & -1 & -2 & $\infty$  \\ \hline
0 & 0  & 0  & 1/2  \\
0 & 1/2 & 1  & 1/2   \\
1 & 1/2 & 1  & 3/2  \\
1 & 1   & 2  &  3/2 
\end{tabular}
\end{table}

Positions of singular fibers of the family $\mc{Y}$ coincide with singular points of its Picard-Fuchs operator. This operator has conifold (type C) singularities at $w=-1,-2$ where the corresponding fibers $Y_{-1},Y_{-2}$ are double octics of different types (see \cite{CS23}). For conifold points, the description of a semistable degeneration and the limiting mixed Hodge structure is given in \cite{oczko2025}. In this paper we study degenerations at $w=0,\infty$. At these two points the Picard-Fuchs operator has singularities of type K.

The Galois representation  at a point of type K is expected to be isomorphic to the Galois representation associated to a weight-three modular form and its Tate twist. Equivalently, this is a Galois representation on the transcendental lattice of a singular K3 surface. Counting points (over finite fields) on the singular fibers $X_0$ and $X_\infty$ we get numerical evidence for the above isomorphism.

The degenerate fiber $X_0$ from the family $\mc{X}$ has the following equation.

\[u^2=x^2yzt(x - y + z) (y - z - t) (x - y - t)\]

This double cover of $\PP^3$ is not normal. The normalization is performed by substituting $u' = \frac{u}{x}$. Geometrically it removes the double plane $\Pi=\{x=0\}$ from the branch divisor. Restricting the equation to $\Pi$, we get the following double sextic.

\[u^2 = yzt(-y+z)(y-z-t)(-y-t)\]

This double sextic can be resolved to a singular K3 surface (\cite{Per85} p. 266) which we denote by $K_0$ its associated modular form has LMFDB label 32.3.d.a. \cite{lmfdb} we denote it as $f_{32}$.

The modular form $f_{32}$ is given by the following eta quotient.

\[
f_{32} =\frac{\eta(4z)^5 \eta(8z)^5}{\eta(2z)^2 \eta(16z)^2}=q \prod_{n=1}^\infty (1-q^{2n})^{-2}(1-q^{4n})^5(1-q^{8n})^5(1-q^{16n})^{-2}
\]

In the following table we collect coefficients $a_p$ of the modular form $f_{32}$ and numbers $N_p$ of points of the reduction $X_0$ at $p$ for primes $2 < p \leq 197$.

\begin{table}[H]
\begin{tabular}{ccc|ccc|ccc}
$p$ & $a_p$ & $N_p$ & $p$ & $a_p$ & $N_p$ & $p$ & $a_p$ & $N_p$  \\ \hline
2 & 0  & 16         & 53 & 0   & 154390    & 127& 0 & 2048892 \\
3 & 2 & 38          & 59 & 82   & 205534   & 131& -62 & 2248678 \\
5 & 0 & 166         & 61 & 0   & 234302    & 137& -238 & 2608856 \\
7 & 0  & 372        & 67 & -62   & 301094  & 139& -206 & 2686382 \\
11 & -14 & 1390     & 71 & 0   & 358196    & 149& 0 & 3352054 \\
13 & 0   & 2510     & 73 & -142   & 399672 & 151& 0 & 3443556 \\
17 & 2   & 5456     & 79 & 0   & 493356    & 157& 0 & 3918878 \\
19 & 34   & 6902    & 83 & -158   & 572278    & 163& 322 & 4331078 \\
23 & 0   & 12260    & 89 & 146    & 720488    & 167& 0 & 4658132 \\
29 & 0   & 26014    & 97 & -94 & 931392     & 173& 0 & 5237230 \\
31 & 0   & 29916    & 101& 0 & 1050502   & 179& 34 & 5736022 \\
37 & 0  & 53318   & 103& 0 & 1093140     & 181& 0 & 5994902 \\
41 & -46 & 72248  & 107& 178 & 1225294     & 191& 0 & 6968636 \\  
43 & -14 & 79694  & 109& 0 & 1318574     & 193& 98 & 7263072 \\
47 & 0   & 104012 & 113& 98 & 1468112     & 197& 0 & 7722598 \\
\end{tabular}
\caption{Point counts for $X_0$ and Fourier coefficients of $f_{32}$.}
\end{table}

We observe the following relation for primes $2 < p \leq 197$.

\[
N_p = 1+\left(1-3\left( \frac{-1}{p}\right)\right)p - a_p+\left(1+\left( \frac{-1}{p}\right)\right)p^2+p^3
\]

For $w=\infty$ the appropriate double sextic has the following equation.

\[u^2 = xzt(x+z)(-z-t)(x-t)\]

This double sextic can also be resolved to a singular K3 surface \cite{Per85} which we denote as $K_\infty$ its associated modular form has LMFDB label 16.3.c.a. \cite{lmfdb}, we denote it as $f_{16}$. This modular form is given by the following eta product.

\[
f_{16}=\eta(4z)^6=q\prod_{n=1}^\infty(1-q^{4n})^6
\]

Let $a_p$ denote the coefficients of $f_{16}$ and by $N_p$ the number of points of reduction of $X_\infty$ modulo $p$ we get the following table.

\begin{table}[H]
\begin{tabular}{ccc|ccc|ccc}
$p$ & $a_p$ & $N_p$ & $p$ & $a_p$ & $N_p$ & $p$ & $a_p$ & $N_p$  \\ \hline
2 & 0  & 16         & 53 & 90   & 151544    & 127& 0 & 2064640 \\
3 & 0 & 40          & 59 & 0  & 208920   & 131& 0 & 2265384 \\
5 & -6 & 152         & 61 & -22   & 230664    & 137& 210 & 2589776 \\
7 & 0  & 400        & 67 & 0 & 305320  & 139& 0 & 2705080 \\
11 & 0 & 1464     & 71 & 0   & 363024    & 149& -102 & 3330104 \\
13 & 10   & 2344     & 73 & -110 & 394384 & 151& 0 & 3465904  \\
17 & -30   & 5216     & 79 & 0   & 499360    & 157& 170 & 3894216 \\
19 & 0   & 7240    & 83 &  0  & 578760    & 163& 0 & 4357480 \\
23 & 0   & 12720    & 89 &   -78  & 712880    & 167& 0 & 4685520 \\
29 & 42   & 25160    & 97 & 130 & 921856     & 173& 330 & 5207144 \\
31 & 0   & 30784    & 101& -198 & 1040600   & 179& 0 & 5767560 \\
37 & -70  & 52056  & 103& 0 & 1103440     & 181& -38 & 5962360 \\
41 & 18 & 70544  & 107& 0 & 1236600     & 191& 0 & 7004544 \\  
43 & 0 & 81400  & 109& -182 & 1306984     & 193& -190 & 7226304 \\
47 & 0   & 106080 & 113& -30 & 1455584     & 197& -390 & 7684376 \\
\end{tabular}
\caption{Point counts for $X_\infty$ and Fourier coefficients of $f_{16}$.}
\end{table}

We observe the following relation for primes $2 < p \leq 197$.

\[
N_p=1-\left(\frac{-1}{p}\right)p-a_p+p^2+p^3
\]

Our goal is to construct the Galois representation on the limiting mixed Hodge structure explicitly and to explain rigorously the relationship between this Galois representation and the Galois representation associated to the appropriate modular form.

The construction of the Galois action is realized through a detailed study of a semistable reduction of the family $\mc{X}$ at the points $w=0$ and $w=\infty$. We concentrate on the more complicated case of $w=0$.

\section{Resolution of a Generic Fiber}

The family of double covers of $\PP^3$ introduced in the previous section has the following equation.

\begin{multline}\label{eq1}
\mc{X}:=\{((x:y:z:t:u),(w_0:w_1))\in \PP(1,1,1,1,4)\times\PP^1: \\
xyzt(x - y + z) (y - z - t) (x - y - t) (w_1x + w_0y)=u^2\}
\end{multline}

Our goal is to obtain a precise description of the singularities of the degenerate fiber $X_0$. We do so by carefully studying how each blow-up used to resolve generic fibers affects the geometry of the degenerate fiber.

We shall visualize this process using octic diagrams introduced in \cite{oczko2025}. For the sake of completeness we briefly recall this approach.

For a double cover of an octic arrangement an octic diagram is a row of squares, with squares representing components in the branching divisor. When there is a line with a label $a$ on a square $b$ it means that components $a$ and $b$ in the branching divisor intersect along a curve. By abuse of notation we use the same labels $S_i$ both for components in the plane arrangement as well as for their strict transforms at each step of the resolution. In diagrams we omit the letter "$S$" when labeling intersections.

We first observe that we only need to study blow-ups of the degenerate fiber $X_0$ which involve planes $S_1$ or $S_8$, as other planes remain in the same configuration as in generic fibers. Planes $S_1$ and $S_8$ are two planes which degenerate to a double plane in the degenerate fiber at $w=0$. The configuration of lines cut out by other components on the double plane $S_1=S_8$ is as follows.

\includegraphics[scale=1.2]{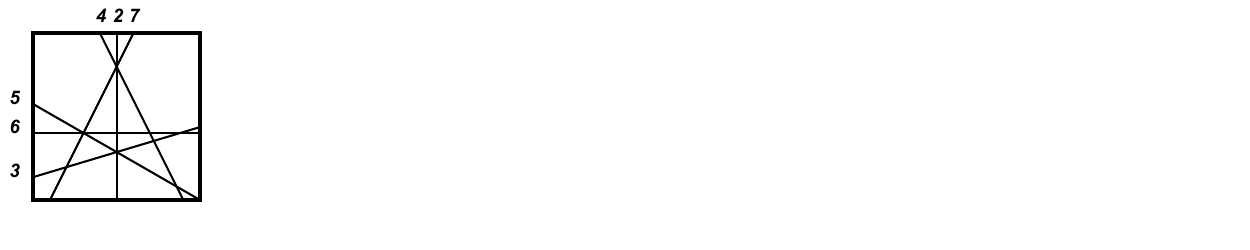}

\textbf{Blow-ups of fivefold points and triple lines.} We begin the resolution by blowing up two $p^1_5$ points: $S_{12358}$ and $S_{12478}$. As the multiplicity of these singularities is odd, these blow-ups introduce new components $A$ and $B$ in the branching divisor. Next we blow up the triple line $S_{128}$ which introduces a new component $C$. After these blow-ups there are no fivefold points or triple lines in generic fibers. In the degenerate fiber the configuration of lines on the double surface is as follows.

\includegraphics[scale=1.2]{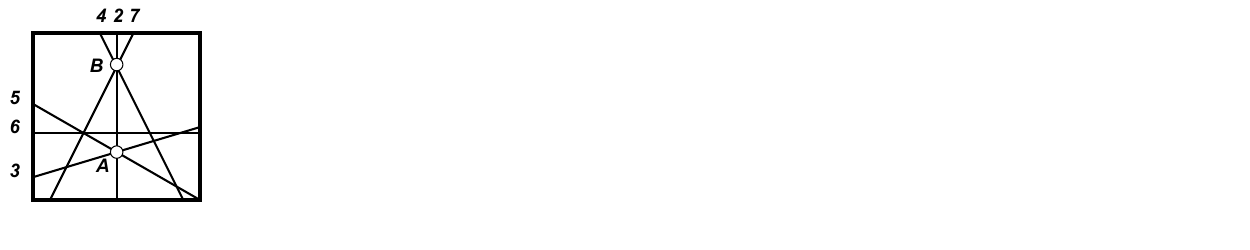}

\textbf{Blow-up of fourfold points.} In the next step, we blow up $p^0_4$ points. The only $p^0_4$ point lying on the double surface is the point $S_{1456}$, which in the degenerate fiber becomes a quintuple point $S_{14568}$. As the multiplicity of this singularity is even in generic fibers, and odd in the degenerate fiber, its blow up introduces a new component only in the degenerate fiber. The center of this blow-up lies on surfaces $S_1,S_4,S_5,S_6$ in generic fibers, while in the degenerate fiber it also lies on the surface $S_8$. As a result the surface $S_8$ splits into two surfaces in the degenerate fiber, one of them being the new component which we denote as $S_{8^0}$ and the other one being the strict transform of the surface $S_8$. Intersections of double planes in the degenerate fiber with other components are as follows.

\includegraphics[scale=1.2]{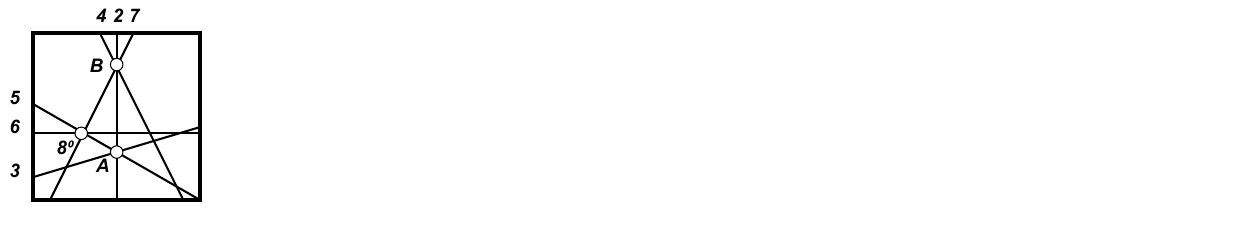}

\textbf{Blow-ups of double lines.} The only singularities in generic fibers left to resolve are double lines. Depending on the order we choose we might obtain different results. Here, we choose the order which leads to the simplest possible calculations. We first blow up all double lines which do not lie on components $S_1$ or $S_8$ in generic fibers. Blow-ups of these double lines do not introduce new singularities and only separate appropriate components of the branch divisor.

We next blow up double lines $S_{14},S_{15},S_{16}$ which become triple lines in the degenerate fiber as they also intersect the surface $S_8$. The surface $S_8$ is not blown up in generic fibers and is blown up in the degenerate fiber, therefore we get new components which appear only in the degenerate fiber. These new components are $S_{8^4},S_{8^5},S_{8^6}$ (which are over lines $S_{14},S_{15}$ and $S_{16}$ respectively). An octic diagram of intersections between components in the degenerate fiber is as follows.

\includegraphics[scale=1.1]{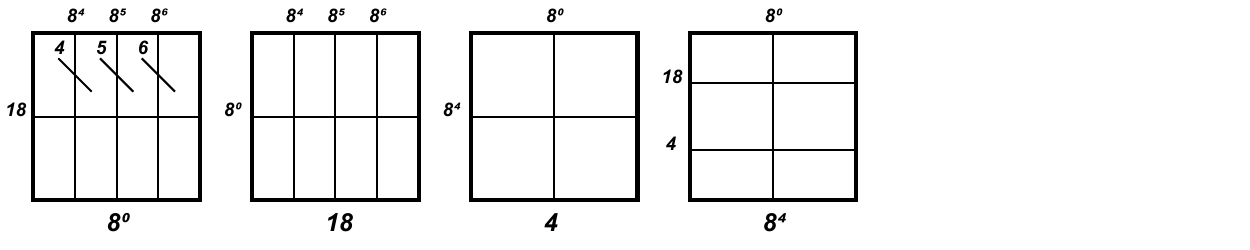}

Next, we blow up double other double lines which lie on the surface $S_1$ in generic fibers, these are $S_{13},S_{17},S_{1A},S_{1B}$ and $S_{1C}$. We introduce new components to the branching divisor as these double lines from generic fibers become triple lines in the degenerate fiber. Let us denote the new components by $S_{8^3},S_{8^7},S_{8^A},S_{8^B},S_{8^C}$ (these are over lines $S_{13},S_{17},S_{1A},S_{1B}$ and $S_{1C}$ respectively).

\includegraphics[scale=1.2]{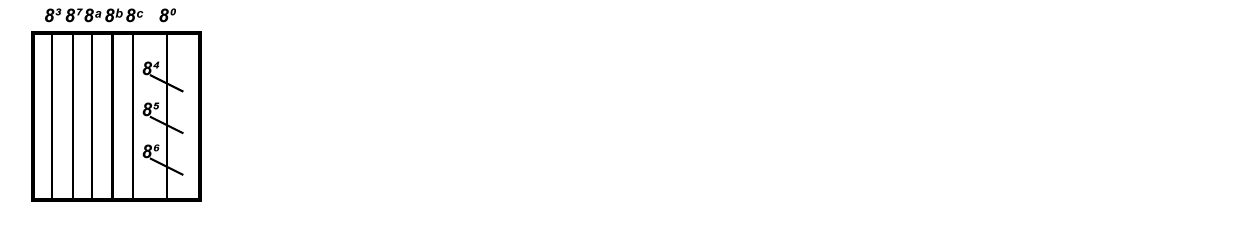}

The only singularities remaining in generic fibers are the double lines lying on $S_8$. In the degenerate fiber these blow-ups separate components $S_{8^0},S_{8^3},S_{8^4},$ $S_{8^5},S_{8^6},S_{8^7},S_{8^A},S_{8^B},S_{8^C}$ from components $S_{1},S_{2},S_{3},S_{4},S_{5},S_{6},S_{7},S_{A},S_{B},S_{C}$. Among the lines blown up are $l_{84},l_{85},l_{86}$. In the degenerate fiber each of them becomes a union of two intersecting lines. Therefore, blowing them up introduces a single pinch point on each of the three double curves $l_{8^0 8^4},l_{8^0 8^5},l_{8^0 8^6}$. We denote the resulting family by $\mc{Y}$.

After all double lines have been blown up, no singularities remain in generic fibers. Therefore, generic fibers of $\mc{Y}$ are smooth Calabi--Yau threefolds. On the other hand, the degenerate fiber $Y_0$ is singular along the double plane $S_1=S_8$ and in three double curves $l_{8^0 8^4},l_{8^0 8^5},l_{8^0 8^6}$. Each of these double curves has a single pinch point and intersects the double plane at a single point.

\section{Semistable reduction}
In the previous section we resolved a generic fiber of $\mc{X}$ and described singularities of the degenerate fiber $Y_0$ of $\mc{Y}$. Now, we shall proceed to construct a semistable degeneration.

\begin{definition}[\mbox{\cite[p.102]{Morr84}}]
A semistable degeneration $\pi \colon \mc{S} \rr \Delta$ is a proper, flat family, such that $\mc{S}$ is a smooth variety, whose generic fibers over $\Delta^*$ are smooth and whose degenerate fiber over $0$ is a union of smooth varieties intersecting transversally along smooth subvarieties.
\end{definition}

Construction of a semistable degeneration consists of two parts. In the first step we compute a log-resolution of $\mc{Y}$ by blowing up singularities of the degenerate fiber $Y_0$. These blow-ups will add non-reduced components to $Y_0$. Then, in the second part we perform the base change to remedy that. To ease the notation we denote the family as $\mc{Y}$ and the degenerate fiber as $Y_0$ at every step of the process.

The degenerate fiber $Y_0$ is singular in the double surface $S_1=S_8$ and in three double curves  $l_{8^0 8^4},l_{8^0 8^5},l_{8^0 8^6}$ with a single pinch point on each of these double curves. To resolve these singularities we perform blow-ups with smooth centers contained in the degenerate fiber. Each such blow-up adds a new component to the degenerate fiber. These components might be non-reduced depending on the multiplicity of the center of a blow-up.

We begin by blowing up the double surface $S_1=S_8$. This double surface is a singularity of the threefold $Y_0$ as well as a singularity of the whole family $\mc{Y}$ viewed as a fourfold. Therefore, when we blow up this surface we resolve a singularity both in the fourfold and in the degenerate fiber. The resulting degenerate fiber is of the form $Y_0 = Y + S$, where $Y$ is the strict transform of $Y_0$ and $S$ is the new component, which is a conic bundle over the blown-up surface. As the fourfold $\mc{Y}$ is now smooth and the threefold $Y$ is normal, the new components arising from blow-ups with centers contained in $Y$ have the same multiplicities as the subvarieties being blown up.

Next, we blow up the double lines $l_{8^0 8^4},l_{8^0 8^5},l_{8^0 8^6}$. This introduces nonreduced components $Q_1,Q_2,Q_3$ of multiplicity two; each is a $\PP^2$-bundle over the corresponding blown-up curve.

The degenerate fiber is now of the form $Y_0 = Y + S +2Q_1+2Q_2+2Q_3$. These components intersect transversely, but are not reduced. Taking a double cover of $\mc{Y}$ branched along $Y\cup S$ would produce a reduced degenerate fiber; however, it would also introduce new singularities because $Y$ and $S$ intersect. To remedy this, we blow up the intersection $Y \cap S$. This adds a new component $V$ to the central fiber which has multiplicity two and separates components $Y$ and $S$.

After these blow-ups the degenerate fiber has the following components.

\begin{itemize}
    \item Strict transform $Y'$ of the threefold $Y$ we started with, it has multiplicity one.
    \item A $\PP^1$ bundle $V$ over a K3 surface $K_0$ with multiplicity two.
    \item A conic bundle $S$ with multiplicity one over the blown up double surface $S_1=S_8$. The surface $S_1=S_8$ is $\PP^2$ blown up at 18 points.
    \item $\PP^2$-bundles $Q_1,Q_2,Q_3$ with multiplicity two which were blown up along a line. These bundles are over double lines $l_{8^0 8^4},l_{8^0 8^5},l_{8^0 8^6}$ respectively.
\end{itemize}

We now take a double cover of the entire family $\mc{Y}$, branched along $S \cup Y$. As a result we get a simple normal crossing divisor in the degenerate fiber and thus we have obtained a semistable family $\mc{Z}$.

In the final family $\mc{Z}$ components $V,Q_1,Q_2$ and $Q_3$ in the degenerate fiber $Z_0$ are replaced by appropriate double covers $V',Q_1',Q_2',Q_3'$. The description of all components in $Z_0$ is the following.

\begin{itemize}
    \item $Y$ is a smooth resolution of a double sextic.
    \item $S$ is a conic bundle over the double surface.
    \item $Q_1',Q_2',Q_3'$ are now quadric bundles.
    \item $V'$ is now a conic bundle over a K3 surface $K_0$. 
\end{itemize}

We have the following nonempty intersections of components of the degenerate fiber.

\begin{itemize}
    \item $Y \cap Q_i'$ are conic bundles over $\PP^1$ each of them having one singular fiber consisting of two intersecting lines.
        \item The surfaces $S \cap Q_i'$ are isomorphic to $\PP^2$.
    \item The surfaces $V' \cap S$ and $V' \cap Y$ are isomorphic to the K3 surface $K_0$.
    \item The surfaces $V' \cap Q_i'$ are isomorphic to a double cover of $\PP^1 \times \PP^1$ branched along two separate lines.
    \item $Y \cap V' \cap Q_i'$ and $S \cap V' \cap Q_i'$ are isomorphic to $\PP^1$.
\end{itemize}

\section{Limiting Mixed Hodge Structure}

\subsection{Monodromy weight spectral sequence}
We compute the limiting mixed Hodge structure of the semistable family $\pi \colon \mc{Z} \rr \Delta$ constructed in the previous section on the degenerate fiber $Z_0$. For that purpose we use the monodromy weight spectral sequence. Its first page consists of cohomology groups of components of the degenerate fiber $Z_0$ and their intersections. It has the following entries:

\[
E^{-k,h+k}_1=\bigoplus_{j\geq \max\{-k,0\}}H^{h-2j-k}(Z_0^{[2j+k+1]},\QQ_\ell)[-j -k]
\]
(for details see \cite[p.273]{Steen08}).

$Z_0^{[n]} 
:= \bigsqcup_{i_1,...,i_n} U_{i_1} \cap ...\cap U_{i_n}$ denotes a disjoint union of $n$-fold intersections of components $U_i$ of the central fiber $Z_0$. The groups $H^{h-2j-k}(Z_0^{[2j+k+1]},\QQ_\ell)$ are Tate twists of $\ell$-adic cohomology groups over $\QQ_\ell$. Horizontal maps come from inclusions and Gysin morphisms. By the following theorem this spectral sequence converges to the limiting mixed Hodge structure $H^q(Z_\infty)$.

\begin{theorem}[{\cite[Cor. 11.23]{Steen08}}]
    The monodromy weight spectral sequence degenerates at the second page $E_2$ and
    \[
    E_1^{-p,q+p} \implies H^q(Z_\infty,\QQ_\ell).
    \]
\end{theorem}

Moreover, on the level of vector spaces the limiting Hodge structure agrees with Hodge structure of generic fibers.

\begin{proposition}[{\cite[Cor. 11.25]{Steen08}}]
    $\operatorname{dim}F^p H^k(Z_\infty) = \operatorname{dim}F^p H^k(Z_t),t 
    \in \Delta^*$
\end{proposition}

In our case there are no quadruple intersections in the degenerate fiber ($Z_0^{[i]} = \varnothing$ for $i \geq 4$), hence the first page of the monodromy weight spectral sequence has the following terms.

\[\begin{tikzcd}[row sep=small]
    H^2(Z_0^{[3]}) & H^4(Z_0^{[2]}) & H^6(Z_0^{[1]}) & 0 & 0\\
    H^1(Z_0^{[3]}) & H^3(Z_0^{[2]}) & H^5(Z_0^{[1]}) & 0 & 0\\
    H^0(Z_0^{[3]}) & H^2(Z_0^{[2]}) & H^4(Z_0^{[1]}) \oplus H^2(Z_0^{[3]}) & H^4(Z_0^{[2]}) & 0\\
    0 & H^1(Z_0^{[2]}) & H^3(Z_0^{[1]}) \oplus H^1(Z_0^{[3]}) & H^3(Z_0^{[2]}) & 0\\
    0 & H^0(Z_0^{[2]}) & H^2(Z_0^{[1]}) \oplus H^0(Z_0^{[3]}) & H^2(Z_0^{[2]}) & H^2(Z_0^{[3]}) \\
    0 & 0 & H^1(Z_0^{[1]}) & H^1(Z_0^{[2]}) & H^1(Z_0^{[3]})\\
    0 & 0 & H^0(Z_0^{[1]}) & H^0(Z_0^{[2]}) & H^0(Z_0^{[3]})\\
        \arrow[from=1-1, to=1-2]
	\arrow[from=1-2, to=1-3]
	\arrow[from=1-3, to=1-4]
	\arrow[from=1-4, to=1-5]
	\arrow[from=2-1, to=2-2]
	\arrow[from=2-2, to=2-3]
	\arrow[from=2-3, to=2-4]
	\arrow[from=2-4, to=2-5]
	\arrow[from=3-1, to=3-2]
	\arrow[from=3-2, to=3-3]
	\arrow[from=3-3, to=3-4]
	\arrow[from=3-4, to=3-5]
	\arrow[from=4-1, to=4-2]
	\arrow[from=4-2, to=4-3]
	\arrow[from=4-3, to=4-4]
	\arrow[from=4-4, to=4-5]
	\arrow[from=5-1, to=5-2]
	\arrow[from=5-2, to=5-3]
	\arrow[from=5-3, to=5-4]
	\arrow[from=5-4, to=5-5]
	\arrow[from=6-1, to=6-2]
	\arrow[from=6-2, to=6-3]
	\arrow[from=6-3, to=6-4]
	\arrow[from=6-4, to=6-5]
	\arrow[from=7-1, to=7-2]
	\arrow[from=7-2, to=7-3]
	\arrow[from=7-3, to=7-4]
	\arrow[from=7-4, to=7-5]
\end{tikzcd}\]  

While the geometric description of components $S,V'$ and $Q_i'$ is relatively easy, we need to consider the threefold $Y$ in more detail. Recall that $Y$ is a resolution of singularities of the double cover of $\PP^3$ branched along the following surface.

\[
x^2yzt(x - y + z) (y - z - t) (x - y - t) =0
\]

After the first stage of the resolution, in which we resolved generic fibers, we obtain a double cover $Y'$ of $\PP^3$ blown up at 6 points and 39 lines. Therefore, the Euler characteristic of the base space of $Y'$ is $4 + 6 \cdot 2 + 39 \cdot 2 = 94$. Its branching locus, excluding the double plane $x=0$, consists of nine copies of $\PP^2$ and nine copies of $\PP^1 \times \PP^1$, which were in total blown up at forty-four points. These components intersect in three triple points, each of them lying on three double lines. Therefore, the Euler characteristic of the branching locus of $Y'$ is $9 \cdot 3 + 9 \cdot 4 + 44 -3 (3 \cdot 2 - 1) = 92$ and the Euler characteristic of $Y'$ is $2 \cdot 94 - 92 = 96$.

The blow-up of a union of two intersecting double lines in $Y'$ increases its Euler characteristic by $3$. The blow-up of (the strict transform of) the double line which passed through the intersection point also increases the Euler characteristic by $3$. Doing these blow-ups for three triple points increases the Euler characteristic by $18$ yielding $\chi(Y) = 18+\chi(Y')=18+96=114$.

Since the blow-ups of $Y$ are centered in rational subvarieties located in the branching divisor, they do not affect its third Betti number, therefore the Betti numbers of $Y$ are $(1,0,56,0,56,0,1)$. Betti numbers of other components and their intersections are as follows.

\begin{itemize}
    \item $V'$ is a conic bundle over a K3 surface, by the Leray spectral sequence $b(V')=(1,0,23,0,23,0,1)$.
    \item Conic bundles $Q_i'$ are double covers of line bundles $Q_i$ branched along a disjoint union of $Q_i \cap S$ and $Q_i \cap Y$. Consequently, $\chi(Q_i')= 2\chi(Q_i) - \chi(Q_i \cap S) - \chi(Q_i \cap Y)=16 -3 -5 =8$, therefore $b(Q_i)=(1,0,3,0,3,0,1)$.
    \item $S$ is the blow-up at three points of a conic bundle over the projective plane $\PP^2$ blown up in 15 points, hence $b(S) = (1,0,20,0,20,0,1)$.
    \item As observed in the previous section $Y \cap Q_i'$ are conic bundles over $\PP^1$ each of them having one singular fiber consisting of two intersecting lines. Therefore $\chi(Y\cap Q_i')=2 \chi(\PP^1 \times \PP^1) - 3=5$. Consequently $b(Y \cap Q_i') =(1,0,3,0,1)$.
    \item $Y \cap V'$ and $S \cap V'$ are isomorphic to the K3 surface $K_0$, therefore $b(Y \cap V') = b(S \cap V') = (1, 0, 22, 0, 1)$.
    \item $S \cap Q_i'$ are isomorphic to $\PP^2$, thus $b(S \cap Q_i') = (1,0,1,0,1)$.
    \item $V \cap Q_i'$ are isomorphic to $\PP^1 \times \PP^1$, thus $b(V \cap Q_i') = (1,0,2,0,1)$.
\end{itemize}

The only nonempty triple intersections are $Y \cap V' \cap Q_i'$ and $Y \cap S' \cap Q_i'$, all six of them are isomorphic to $\PP^1$. Putting this data into the monodromy weight spectral sequence we get the following first page.

\[\begin{tikzcd}[row sep=small]
    \CC^6 & \CC^{11} & \CC^6 & 0 & 0\\
    0 & 0 & 0 & 0 & 0\\
    \CC^6 & \CC^{62} & \CC^{108} \oplus \CC^6 & \CC^{11} & 0\\
    0 & 0 & 0 & 0 & 0\\
    0 & \CC^{11} & \CC^{108} \oplus \CC^6 & \CC^{62}& \CC^6 \\
    0 & 0 & 0 & 0 & 0\\
    0 & 0 & \CC^6 & \CC^{11} & \CC^6\\
        \arrow[from=1-1, to=1-2]
	\arrow[from=1-2, to=1-3]
	\arrow[from=1-3, to=1-4]
	\arrow[from=1-4, to=1-5]
	\arrow[from=2-1, to=2-2]
	\arrow[from=2-2, to=2-3]
	\arrow[from=2-3, to=2-4]
	\arrow[from=2-4, to=2-5]
	\arrow[from=3-1, to=3-2]
	\arrow[from=3-2, to=3-3]
	\arrow[from=3-3, to=3-4]
	\arrow[from=3-4, to=3-5]
	\arrow[from=4-1, to=4-2]
	\arrow[from=4-2, to=4-3]
	\arrow[from=4-3, to=4-4]
	\arrow[from=4-4, to=4-5]
	\arrow[from=5-1, to=5-2]
	\arrow[from=5-2, to=5-3]
	\arrow[from=5-3, to=5-4]
	\arrow[from=5-4, to=5-5]
	\arrow[from=6-1, to=6-2]
	\arrow[from=6-2, to=6-3]
	\arrow[from=6-3, to=6-4]
	\arrow[from=6-4, to=6-5]
	\arrow[from=7-1, to=7-2]
	\arrow[from=7-2, to=7-3]
	\arrow[from=7-3, to=7-4]
	\arrow[from=7-4, to=7-5]
\end{tikzcd}\]

\subsection{Computation of the second page.}
The monodromy weight spectral sequence degenerates at the second page and converges to the cohomology of generic fibers. In our case generic fibers $Y_t$ are Calabi--Yau threefolds with $b(Y_t)=(1,0,49,4,49,0,1)$, (Arr No.70 in \cite{CYK}). We are interested in the Galois representation on the limiting mixed Hodge structure on the third cohomology. Since $\dim H^3_{\mathrm{lim}}(Y_t)=\dim H^3(Y_t)=4$ and the terms on the second page in even rows $E_2^{2,1},E_2^{0,3},E_2^{-2,5}$ vanish, we have that $E_2^{1,2}=E_2^{-1,4}=\CC^2$. By definition, $E_2^{-1,4}$ is the quotient 
\[\operatorname{ker}(E^1_{-1,4} \rr E^1_{0,4})/\operatorname{Im}(E^1_{-2,4}\rr E^1_{-1,4}).
\]
To compute the Galois representation we need to find two cocycles in the kernel of $\phi:E_1^{-1,4} = \CC^{62} \rr \CC^{108} \oplus \CC^6=E_1^{0,4}$ which are not in the image of $H^0(Z^{[3]}) \rr H^2(Z^{[2]})$ and study their Galois representation. The other term on the second page $E_2^{1,2}$ which contributes to the third cohomologies can be computed by dualizing arguments for $E_2^{-1,4}$.

Recall that $Y \cap V'$ and $S \cap V'$ are isomorphic to the K3 surface $K_0$ which is birational to the double cover of $\PP^2$ branched along the union of six lines $yzt(y-z)(y-z-t)(y+t)=0$. This configuration of six lines has three triple points $(0:0:1),(0:1:0)$ and $(1:1:0)$, thus the even part of the Néron--Severi groups has rank 19. The line $z+t=0$ intersects the union of 6 lines at three double points. Consequently, the pullback of the line $z+t=0$ to the double cover splits into a union of two lines. The difference of these two lines is an odd cycle and the rank of the Picard group is 20. Hence, transcendental lattices of $Y \cap V'$ and $S \cap V'$ generate a four-dimensional vector subspace $M$ of $H^2(Z^{[2]})$. We claim that the image of $M$ under $\phi$ is two-dimensional.

We first observe that all cohomology classes in $H^2(S)$ and $H^2(Y)$ are algebraic. In the case of $S$ we have 18 classes from blow-ups and two algebraic classes on a product of $\PP^2$ and a conic. These are linearly independent in $H^2(S)$ and therefore generate it. For $Y$ we observe that after resolving generic fibers it was a double cover of $\PP^3$ blown up in 39 lines and 6 points. Together with the pullback of a hyperplane class on $\PP^3$, this gives 46 linearly independent classes. That threefold was then blown up along nine lines, which introduces nine new algebraic classes independent from the rest. The cohomology group $H^2(Y)$ has dimension 56, and we have found 55 independent elements in its subspace $H^{1,1}(Y)$, i.e. $h^{1,1}(Y) \geq 55$. By Hodge theory, $H^2(Y) \cong H^{2,0}(Y) \oplus H^{1,1}(Y) \oplus H^{0,2}(Y)$, and $H^{2,0}(Y) \cong H^{0,2}(Y)$, thus $h^{1,1}(Y) =56$. For $V'$ the Künneth formula gives that $H^2(V')$ has dimension 23 and its subgroup generated by algebraic cycles has dimension 21, hence the transcendental part of $H^2(V')$ is two-dimensional.

It follows that we have demonstrated that the transcendental subspace of $H^2(Z^{[2]})$ is generated by transcendental lattices on two K3 surfaces and is four-dimensional. On the other hand we have shown that transcendental subspace of $H^4(Z^{[1]})$ is two-dimensional. Hence, the kernel of $\phi$ consists of transcendental cycles on the K3 surface $K_0$ and is two-dimensional. Consequently, the second page of the monodromy weight spectral sequence is as follows.

\[\begin{tikzcd}[row sep=small]
    0 & 0 & \CC & 0 & 0\\
    0 & 0 & 0 & 0 & 0\\
    0 & \CC^2 & \CC^{49}& 0 & 0\\
    0 & 0 & 0 & 0 & 0\\
    0 & 0 & \CC^{49} & \CC^2& 0\\
    0 & 0 & 0 & 0 & 0\\
    0 & 0 & \CC & 0 & 0\\
        \arrow[from=1-1, to=1-2]
	\arrow[from=1-2, to=1-3]
	\arrow[from=1-3, to=1-4]
	\arrow[from=1-4, to=1-5]
	\arrow[from=2-1, to=2-2]
	\arrow[from=2-2, to=2-3]
	\arrow[from=2-3, to=2-4]
	\arrow[from=2-4, to=2-5]
	\arrow[from=3-1, to=3-2]
	\arrow[from=3-2, to=3-3]
	\arrow[from=3-3, to=3-4]
	\arrow[from=3-4, to=3-5]
	\arrow[from=4-1, to=4-2]
	\arrow[from=4-2, to=4-3]
	\arrow[from=4-3, to=4-4]
	\arrow[from=4-4, to=4-5]
	\arrow[from=5-1, to=5-2]
	\arrow[from=5-2, to=5-3]
	\arrow[from=5-3, to=5-4]
	\arrow[from=5-4, to=5-5]
	\arrow[from=6-1, to=6-2]
	\arrow[from=6-2, to=6-3]
	\arrow[from=6-3, to=6-4]
	\arrow[from=6-4, to=6-5]
	\arrow[from=7-1, to=7-2]
	\arrow[from=7-2, to=7-3]
	\arrow[from=7-3, to=7-4]
	\arrow[from=7-4, to=7-5]
\end{tikzcd}\]

Consequently, the weight filtration on $H^3(Z_\infty)$ has two nontrivial graded pieces $\operatorname{Gr}_W^2$ and $\operatorname{Gr}_W^4$ and $H^3(Z_\infty)=\operatorname{Gr}_W^2\oplus\operatorname{Gr}_W^4$, up to semisimplification. The explicit description of the monodromy weight spectral sequence implies the following isomorphisms $\operatorname{Gr}_W^4\simeq T(K_0)(-1)$, $\operatorname{Gr}_W^2\simeq T(K_0)$ where $T(K_0)$ denotes the transcendental lattice of the singular K3 surface $K_0$. The vector spaces $\operatorname{Gr}_W^2$ and $\operatorname{Gr}_W^4$ are naturally endowed with Galois actions induced by the Galois action on components of the central fiber.

Since singular K3 surfaces are modular \cite{Liv}, we have a modular form $g$ associated with $K_0$. As $Gr_W^4$ is a Tate twist of $Gr_W^2$ (\cite{Liv}, Example 1.6), we obtain the following identities of $L$-functions:

\[
L(Gr_W^2,s)=L(g,s), \quad L(Gr_W^4,s) = L(g,s-1)
\]
It follows that 

\[L(H^3_{lim}(Z),s)=L(g,s)L(g,s-1).\]

Finally for good odd primes the Frobenius polynomials of two-dimensional pieces are
$1-a_pT+\chi_{-8}(p)p^2T^2$ and $1-a_pT+\chi_{-4}(p)p^2T^2$ for the forms $f_{32}$ and $f_{16}$, respectively. After the Tate twist they become
$1-pa_pT+\chi_{-8}(p)p^4T^2$ and $1-pa_pT+\chi_{-4}(p)p^4T^2$, respectively. These two-dimensional pieces come precisely from transcendental cycles on K3 surfaces. In particular, the Frobenius polynomial of $H^3_{\mathrm{lim}}(Z)$ is the product of these two quadratic factors.

The main theorem now follows from preceding computations.

\begin{theorem}\label{main thm}
    Let $\mc{Y}$ be a family of Calabi-Yau threefolds (family No. 70 from \cite{CYK}), with two singular fibers $Y_0,Y_\infty$ of type K. Let $f_1,f_2$ be two modular forms 32.3.d.a and 16.3.c.a. of weight-three associated to K3 surfaces \cite{lmfdb}. The Galois representation on the limiting mixed Hodge structure of $Y_0$ (resp. $Y_\infty$) is isomorphic up to semisimplification to the direct sum of the Galois representation associated to $f_1$ (resp. $f_2$) and its Tate twist.
\end{theorem}

\printbibliography

@article{CYK,
author = {Slawomir Cynk and Beata Kocel-Cynk},
year = {2019},
title = {Classification of Double Octic Calabi-Yau Threefolds with $h^{1,2} \leq 1$ defined by an Arrangement of Eight Planes},
volume = {22},
pages={1--52},
journal = {Communications in Contemporary Mathematics}
}

@article{oczko2025,
author = {Oczko, Marcin},
year = {2026},
title = {Semistable degenerations of double octics},
volume = {177},
journal = {manuscripta mathematica}
}

@misc{lmfdb,
  author       = {The LMFDB Collaboration},
  title        = {The L-functions and modular forms database},
  url = {https://www.lmfdb.org},
  year         = {2025},
  note         = {[Online; accessed 27 November 2025]},
}

@ARTICLE{CO24,
  title     = "Counterexample to the {N{\'e}ron--Ogg--Shafarevich} criterion
               for {Calabi--Yau} threefolds",
  author    = "Chmiel, T. and Oczko, M.",
  journal   = "Journal of Algebra",
  volume    =  "701",
  pages     = "110--123",
  year      =  2026
}

@book{Steen08,
author = {Chris A.M. Peters and Joseph H.M. Steenbrink},
year = {2008},
pages = {},
publisher = "Springer",
title = {Mixed Hodge Structures}
}

@inbook{Morr84,
title = {Chapter VI. The Clemens-Schmid exact sequence and applications},
booktitle = {Topics in Transcendental Algebraic Geometry. (AM-106), Volume 106},
author = {David R. Morrison},
publisher = {Princeton University Press},
pages = {101--120},
year = {1984}
}

@article{PP81,
  title={Degeneration of surfaces with trivial canonical bundle},
  author={Ulf Persson and Henry C. Pinkham},
  journal={Annals of Mathematics},
  year={1981},
  volume={113},
  pages={45--66}
}

@article{CvS21,
author = {Sławomir Cynk and Duco van Straten},
title = {A special Calabi–Yau degeneration with trivial monodromy},
journal = {Communications in Contemporary Mathematics},
year = {2022},
volume = {24},
pages={2150055}
}

@article{CS23,
author = {Chmiel, Tymoteusz and Cynk, Slawomir},
year = {2023},
title = {Periods of singular double octic Calabi–Yau threefolds and modular forms},
volume = {296},
journal = {Mathematische Nachrichten},
pages ={3257–3271}
}

@InProceedings{Per85,
author="Persson, Ulf",
title="Double sextics and singular K-3 surfaces",
booktitle="Algebraic Geometry Sitges (Barcelona) 1983",
year="1985",
publisher="Springer Berlin Heidelberg",
pages="262--328",
}

@article{Liv,
author = {Livné, Ron},
year = {1995},
pages = {149-156},
title = {Motivic orthogonal two-dimensional representations of Gal($\wb{\QQ}/\QQ$)},
volume = {92},
journal = {Israel Journal of Mathematics}
}

@article{Orph,
author = {Cynk, Slawomir and Straten, Duco},
year = {2017},
title = {Picard-Fuchs operators for octic arrangements I (The case of orphans)},
volume = {13},
journal = {Communications in Number Theory and Physics}
}

@book{SGA7,
  title={Groupes de Monodromie en Géométrie Algébrique: SGA 7},
  author={Grothendieck, A. and Deligne, P. and Katz, N.},
  year={1972},
  publisher={Springer-Verlag}
}

@article{Landman,
author = {Landman, Alan},
year = {1973},
pages = {89-126},
title = {On the Picard-Lefschetz Transformation for Algebraic Manifolds Acquiring General Singularities},
volume = {181},
journal = {Transactions of The American Mathematical Society - TRANS AMER MATH SOC},
}

@article{Kul77,
author = {Viktor Kulikov},
year = {1977},
pages = {957},
title = {Degenerations of K3 surfaces and Enriques surfaces},
volume = {11},
journal = {Mathematics of the USSR-Izvestiya}
}

@article{Kul81,
author = {Viktor Kulikov},
year = {1981},
pages = {339},
title = {On modifications of degenerations of surfaces with $\kappa = 0$},
volume = {17},
journal = {Mathematics of the USSR-Izvestiya}
}

@article{CLL,
author = {Chiarellotto, Bruno and Lazda, Christopher and Liedtke, Christian},
year = {2019},
pages = {469-514},
title = {A Néron–Ogg–Shafarevich criterion for K3 surfaces},
volume = {119},
journal = {Proceedings of the London Mathematical Society}
}
\end{document}